\documentclass[a4paper]{amsart}
\usepackage{amssymb}
\input xy
\xyoption{all}

\newcommand{\Id}{\mathrm{Id}}
\newcommand{\Comp}{\mathsf{Comp}}

\newcommand{\id}{\mathrm{id}}

\newcommand{\I}{\mathbb I}
\newcommand{\R}{\mathbb R}
\newcommand{\U}{\mathbb U}
\newcommand{\OM}{\mathbb O}
\newcommand{\N}{\mathbb N}

\newcommand{\F}{\mathcal F}
\newcommand{\C}{\mathcal C}

\newcommand{\E}{\mathbb E}
\newcommand{\bul}{{\boldsymbol{\cdot}}}

\newtheorem{theorem}{Theorem}
\newtheorem{df}{Definition}
\newtheorem{lemma}{Lemma}

\begin{document}

\title{On the functor of comonotonically maxitive functionals}

%(revisited version 1 (January 2023))

\author{Taras Radul}

\maketitle

Institute of Mathematics, Kazimierz Wielki University in Bydgoszcz, Poland;
\newline
Department of Mechanics and Mathematics, Ivan Franko National University of Lviv,
Universytettska st., 1. 79000 Lviv, Ukraine.
\newline
e-mail: tarasradul@yahoo.co.uk



\begin{abstract} We introduce a functor of functionals which preserve maximum of comonotone functions and addition of constants. This functor is a subfunctor of the functor of order-preserving functionals and contains  the idempotent measure functor  as subfunctor. The main aim of this paper is to show that this functor is isomorphic to the capacity functor. We establish such isomorphism using the fuzzy max-plus integral.   In fact,  we can consider this result as an idempotent  analogue of  Riesz Theorem about a correspondence between the set of $\sigma$-additive regular Borel measures and the set of linear  positively defined functionals.
 \end{abstract}

\maketitle

\section{Introduction} The general theory of functors acting in the
category $\Comp$ of compact Haussdorff spaces (compacta) and continuous
mappings was founded by E.V.~Shchepin \cite{Shchepin}. He distinguished some elementary
properties of such functors and defined the notion of normal functor that
has become very fruitful.
The class of normal functors and close to them
includes many classical topological constructions: hyperspace $\exp$ , space of
probability measures $P$, superextension $\lambda$, space of hyperspaces of
inclusion $G$ and many other functors (see for example the review \cite{FZ}).

Some functors in the category $\Comp$ categories of topological spaces and continuous maps have also natural algebraical structures. Such structures can be described by the notion of  monad (or triple) structure in the
sense of S.Eilenberg and J.Moore  and their corresponding category of Eilenberg-Moore algebras \cite{EilMoore}.
Many classical functors can be completed to monads: hyperspace \cite{Wyl}, space
of probability measures \cite{Swir}, superextension  \cite{Za}, hyperspaces of
inclusion \cite{Ra}, idempotent measures \cite{Za1}, capacities \cite{NZ}, order-preserving functionals \cite{Ra1} etc.

It seems that the main problem to obtain general results in the theory of monads in the category $\Comp$ is the different nature of functors. There were introduced in \cite{Ra2} and \cite{Ra3} sufficiently wide classes of monads which  have a functional representation, i.e., their functorial part $FX$ can be naturally with preserving of monad structures imbedded in $\R^{CX}$.  Functorial parts of some of such monads are by the definition subspaces of $\R^{CX}$. The most known are the probability measure monad and its idempotent analogue the idempotent (max-plus) measure monad.  Let us remark that  the term 'measure' is used  for  functionals in both cases. In the case of probability measures it is justified by the Riesz Theorem which establishes  a correspondence between the set of $\sigma$-additive regular Borel measures on compacta  and the set of linear  positively defined functionals. This correspondence, facilitated by the Lebesgue integral, allows the term 'measure' to be used for these functionals as well. (See, for example the review   \cite{Fed2} devoted to probability measures on topological spaces where measures generally are considered as functionals.) Recently, a correspondence was obtained in \cite{Ra4} and \cite{Ra5} between probability capacities and functionals preserving maximum of functions and addition of constants referred in \cite{Za1} as 'idempotent measure'.   This correspondence, facilitated by the introduced in \cite{Ra5} max-plus integral, justifies the application of the term 'measure' to these functionals and can be seen as an idempotent counterpart of the Riesz Theorem.

Capacities (non-additive measures, fuzzy measures) were introduced by Choquet in \cite{Ch} as a natural generalization of additive measures. They found numerous applications (see for example \cite{EK},\cite{Gil},\cite{Sch}, \cite{Grab}).  Additionally, Zadeh introduced possibility measures (normalized capacities that preserve the maximum) and founded possibility theory \cite{Zad}, which has since been extensively developed and found numerous applications, see for example \cite{DP},\cite{DP1}. Capacities on compacta were considered in \cite{Lin} where the important role plays the upper-semicontinuity property which  connects the capacity theory with the topological structure. Categorical and topological properties of spaces of upper-semicontinuous normalized capacities on compact Hausdorff spaces were investigated in \cite{NZ}.

In fact, the most of applications of non-additive measures to game theory, decision making theory, economics etc deal not with measures as set functions  but with integrals which allow to obtain expected utility or expected pay-off.  Several types of integrals with respect to non-additive measures were developed for different purposes (see for example books \cite{Grab} and  \cite{Den}). Such integrals are called fuzzy integrals. The most known are the Choquet integral based on the addition and the multiplication operations \cite{Ch} and the Sugeno integral  based on the maximum and the minimum operations \cite{Su}. The max-plus integral based on the maximum and the addition operations was introduced in \cite{Ra5} and plays a crucial role in establishing of an isomorphism between possibility capacity monad and idempotent measure monad.

It was remarked in \cite{Ra5} that  the max-plus  integral can be used not only for possibility capacities, but fol all capacities.As a specific outcome of that fact, we obtain a broader space of functionals than the space of max-plus idempotent measures. The problem of  exploring the topological and categorical properties of this structure was posed in \cite{Ra5}.

The monad of order-preserving functionals which contain both probability measure and and idempotent measure monads as submonads was introduced in \cite{Ra1}.

\section {Preliminaries and definitions}\label{Pr}
In what follows, all spaces are assumed to be compacta (compact Hausdorff space) except for $\R$ and the spaces of continuous functions on a compactum. All maps are assumed to be continuous.

We shall denote the
Banach space of continuous functions on a compactum  $X$ endowed with the sup-norm by $C(X)$. For any $\gamma\in\R$ we shall denote the
constant function on $X$ taking the value $\gamma$ by $\gamma_X$. We also consider the natural lattice operations $\vee$ and $\wedge$  on $C(X)$.

First we recall some properties of functionals. Mainly, we follow the terminology from \cite{Grab}. Let $X$ be a compactum.  We call two functions $\varphi$, $\psi\in C(X)$ comonotonic if $(\varphi(x_1)-\varphi(x_2))\cdot(\psi(x_1)-\psi(x_2))\ge 0$ for each $x_1$, $x_2\in X$. Let us remark that a constant function is comonotonic to any function $\psi\in C(X)$.

Let $X$ be a compact space. We shall say that a functional $\mu:C(X)\to\R$
\begin{itemize}
\item {\em is normalized} if $\mu(\gamma_X)=\gamma$ for any $\gamma\in\R$;
\item {\em is monotone} if $\mu(f)\le\mu(g)$ for any functions $f\le g$ in $C(X)$;
\item {\em is maximative}  if $\mu(f\vee g\})=\mu(f)\vee\mu(g)$ for any functions $f,g\in C(X)$;
\item  {\em is plus-homogeneuos} if $\mu(f+\gamma_X)=\mu(f)+\gamma$ for any $f\in C(X)$ and   $\gamma\in \R$;
 \item {\em is comonotonically maximative}  if $\mu(f\vee g)=\mu(f)\vee\mu(g)$ for any comonotonic functions $f,g\in C(X)$
\end{itemize}

Evidently, each maximative functional is monotone. The problem is not so plain for comonotonically maximative functionals. It is known that each comonotonically maximative functional is monotone for finite compacta. The implication was proved for any compactum with some additional conditions on functional. But generally the problem is still open (see \cite{Ra6} for more details).

Normalized monotone and  plus-homogeneuos functionals were referred as order-preserving functionals and the set of all order-preserving functionals on $C(X)$ for a compactum $X$ was denoted by $OX$ in \cite{Ra1}. Normalized maximative and  plus-homogeneuos functionals were referred as idempotent measures and the set of all idempotent measures  for a compactum $X$ was denoted by $IX$ in \cite{Za1}. We have $IX\subset OX$ for each compactum $X$. Both sets $IX$ and $OX$ we consider as subspaces in $\R^C(X)$. It was shown in \cite{Ra1} and \cite{Za1} the the spaces $IX$ and $OX$ are compacta.

We need the definition of capacity on a compactum $X$. We follow a terminology of \cite{NZ}. By $\F(X)$ we denote the family of all closed subsets of a compactum $X$.

\begin{df}\cite{NZ} A function $c:\F(X)\to [0,1]$  is called an {\it upper-semicontinuous capacity} on $X$ if the three following properties hold for each closed subsets $F$ and $G$ of $X$:

1. $c(X)=1$, $c(\emptyset)=0$,

2. if $F\subset G$, then $c(F)\le c(G)$,

3. if $c(F)<a$ for $a\in[0,1]$, then there exists an open set $O\supset F$ such that $c(B)<a$ for each compactum $B\subset O$.
\end{df}

By $MX$ we denote the set  of all upper-semicontinuous  capacities on a compactum $X$. Since all capacities we consider here are upper-semicontinuous, in the following we call elements of the set $MX$ simply capacities.

It was proved in \cite{NZ} that the space $MX$ of all upper-semicontinuous  capacities on a compactum $X$ is a compactum as well, if a topology on $MX$ is defined by a subbase that consists of all sets of the form $O_-(F,a)=\{c\in MX\mid c(F)<a\}$, where $F$ is a closed subset of $X$, $a\in [0,1]$, and $O_+(U,a)=\{c\in MX\mid c(U)=\sup\{c(K)\mid K$ is a compact subset of $U\}>a\}$, where $U$ is an open subset of $X$, $a\in [0,1]$. Since all capacities we consider here are upper-semicontinuous, in the following we call elements of $MX$ simply capacities.

A capacity $c\in MX$ for a compactum $X$ is called  a possibility capacity if for each family $\{A_t\}_{t\in T}$ of closed subsets of $X$ such that $\bigcup_{t\in T}A_t$ is a closed subset of $X$ we have  $c(\bigcup_{t\in T}A_t)=\sup_{t\in T}c(A_t).$
(See \cite{WK} for more details.) We denote by $\Pi X$ a subspace of $MX$ consisting of all possibility capacities. Since $X$ is compact and $c$ is upper-semicontinuous, $c\in \Pi X$ iff $c$ satisfies the simpler requirement that $c(A\cup B)=\max\{c(A),c(B)\}$ for each closed subsets $A$ and $B$ of $X$.  It is easy to check that  $\Pi X$ is a closed subset of $MX$.

Now we consider one of fuzzy integrals, namely the max-plus integral. We denote $\varphi_t=\{x\in X\mid \varphi(x)\ge t\}$ for $\varphi\in C(X)$ and $t\in\R$. We also put $\ln(0)=-\infty$ and $=-\infty+\gamma=-\infty$ for each $\gamma\in\R$.

\begin{df}\cite{Ra5} Let $\varphi\in  C(X)$ be a function and $c\in MX$. The max-plus integral of $\varphi$ w.r.t. $c$ is given by  the formula $$\int_X^{\vee+} \varphi dc=\max\{\ln(c(\varphi_t))+ t\mid t\in\R\}.$$
\end{df}

The following characterization of the max-plus integral with respect possibility capacities  was obtained in \cite{Ra5}: $\nu\in IX$ iff there exists a unique possibility capacity $c\in\Pi X$ such that $\nu(\varphi)=\int_X^{\vee+} \varphi dc$ for any $\varphi\in C(X)$.

\section {A characterization of the max-plus integral}
The main aim of this section is to obtain a characterization of the max-plus integral with respect all capacities.

The following lemma follows from Lemma 4.27 in \cite{Grab} but it also could be proved by  easy  checking.

\begin{lemma}\label{Comon} Let  $\varphi$, $\psi\in C(X)$ be two comonotonic functions. Then we have $\varphi_t\subset\psi_t$ or $\varphi_t\supset\psi_t$ for each $t\in\R$.
\end{lemma}

Consider any capacity $c\in MX$. Let  $I$ be the max-plus integral with respect
to $c$, i.e. $I:C(X)\to\R$ is a functional defined by the formula  $I(\varphi)=\int_X^{\vee+} \varphi dc$ for $\varphi\in C(X)$.

\begin{lemma}\label{Bound} We have $\min_{x\in X}\varphi(x)\le I(\varphi)\le\max_{x\in X}\varphi(x)$ for any $\varphi\in C(X)$.
\end{lemma}

\begin{proof} Put $a=\min_{x\in X}\varphi(x)$ and  $b=\max_{x\in X}\varphi(x)$. The we have  $$a=0+a=c(\varphi_a)+ a\le\max\{\ln(c(\varphi_t))+ t\mid t\in\R\}=$$
$$=\max\{\ln(c(\varphi_t))+ t\mid t\in(-\infty,b]\}\le0+b=b.$$
\end{proof}

\begin{lemma}\label{prop}  The functional $I$ is normalized, monotone, comonotonically maximative and plus-homogeneous.
\end{lemma}

\begin{proof} Lemma \ref{Bound} yields that $I$ is normalized.

Consider any functions $\varphi$, $\psi\in C(X)$ such that $\varphi\le\psi$. The inequality $I(\varphi)\le I(\psi)$ follows from the obvious inclusion $\varphi_t\subset\psi_t$ and monotonicity of $c$.

Let  $\varphi$, $\psi\in C(X)$ be two comonotone functions. The inequality $I(\psi\vee\varphi)\ge I(\psi)\vee I(\varphi)$ follows from the monotonicity of $I$. We have $\nu(\psi_t)+ t\le I(\psi)\vee I(\varphi)$  and $\nu(\varphi_t)+ t\le I(\psi)\vee I(\varphi)$ for each $t\in \R$. Lemma \ref{Comon} yields that $(\psi\vee\varphi)_t=\psi_t$ or $(\psi\vee\varphi)_t=\varphi_t$. Hence $I(\psi\vee\varphi)\le I(\psi)\vee I(\varphi)$ and we proved that $I$ is comonotonically maximative.

Consider any $\alpha\in\R$ and $\psi\in C(X)$. Then we have $$I(\psi+\alpha_X)=\max\{\ln(c((\psi+\alpha_X)_t)+ t\mid t\in\R\}=$$ $$=\max\{\ln(c(\psi_{t-\alpha})+ t-\alpha+\alpha\mid t\in\R\}=\max\{\ln(c(\psi_s)+ s+\alpha\mid s\in\R\}=I(\psi)+\alpha.$$
\end{proof}

For $A\in\F(X)$ put $\Upsilon^t_A=\{\varphi\in C(X,(-\infty,t])\mid \varphi(a)=t$ for each $a\in A\}$. If $A=\emptyset$ we put $\Upsilon^t_A=C(X)$. We denote $\Upsilon^0_A=\Upsilon_A$ for $t=0$. 

Let $\varphi$ and $\psi$ be two function in $C(X)$. We say that $\varphi$ {\it refines} $\psi$ if for each $t\in \R$ there exists $s\in\R$ such that $\varphi^{-1}(t)\subset \psi^{-1}(s)$. 

The following lemma is a slight modification of Lemma 2 from \cite{Ra7}.

\begin{lemma}\label{Comon1} Let  $\varphi\in C(X)$, $t\in\varphi(X)$ and $\psi\in\Upsilon^t_{\varphi_t}$. Then there exists $\psi'\in\Upsilon^t_{\varphi_t}$ such that $\psi'\le\psi$, $\psi'$ refines $\varphi$  and $\psi'$ is comonotone with $\varphi$.
\end{lemma}

Let us remark that there is no condition that $\psi'$ refines $\varphi$  in the formulation of  Lemma 2 in  \cite{Ra7}, but it follows immediately from the construction of $\psi'$ in the proof.

\begin{theorem}\label{repr} A functional $I:C(X)\to\R$ is  normalized, monotone,  comonotonically maximative and plus-homogeneous
if and only if there exists a unique capacity $c\in MX$ such that   $I(\varphi)=\int_X^{\vee+} \varphi dc$ for each $\varphi\in C(X)$.
\end{theorem}

\begin{proof} Sufficiency is proved in Lemma \ref{prop}.

Necessity. Define the function $c:\F(X)\to \R$ by the formula $c(A)=\inf\{e^{I(\varphi)}\mid \varphi\in \Upsilon_A\}$. It is easy to check that $c$ satisfies Conditions 1 and 2 from the definition of capacity.

Let $c(A)<\eta$ for some $\eta\in [0,1]$ and $A\in\F(X)$. Then there exists $\varphi\in \Upsilon_A$ such that $I(\varphi)<\eta$. Put $\varepsilon=\eta-I(\varphi)$. Put $V=\{x\in X\mid \varphi(x)>-\varepsilon/2\}$. Evidently $V$ is an open set containing $A$ as subset. Take any compactum $B\subset V$ and put $\psi=\min\{0_X,\varphi+\varepsilon/2_X\}$. Evidently, $\psi\in \Upsilon_B$. Since $I(\varphi+\varepsilon/2_X)=I(\varphi)+\varepsilon/2$ and $\psi$ is comonotonic with $\varphi+\varepsilon/2_X$ with $\psi\le\varphi+\varepsilon/2_X$, we have $I(\psi)\le I(\varphi)+\varepsilon/2$. Thus, $c(B)\le I(\psi)\le I(\varphi)+\varepsilon/2<\eta$. Hence $c$ is upper semi continuous.

Let us show that $\int_X^{\vee+} \varphi dc=I(\varphi)$ for each $\varphi\in C(X)$.  We have 
$$\int_X^{\vee+} \varphi dc=\max\{\ln(\inf\{e^{I(\chi)})\mid \chi\in \Upsilon_{\varphi_t}\}+t\mid t\in\R\}=$$ 
$$=\max\{\inf\{I(\chi))\mid \chi\in \Upsilon_{\varphi_t}\}+t\mid t\in\R\}=$$
$$=\max\{\inf\{I(\chi+t_X))\mid \chi\in \Upsilon_{\varphi_t}\}\mid t\in\R\}=$$
$$=\max\{\inf\{I(\zeta))\mid \zeta\in \Upsilon^t_{\varphi_t}\}\mid t\in\R\}.$$

Consider the function $\zeta=\min{t_X,\varphi}\in \Upsilon^t_{\varphi_t}$. Since $\zeta$ is comonotonic with $\varphi$ and $\zeta\le\varphi$, we have $I(\zeta)\le I(\varphi)$. Hence $\int_X^{\vee+} \varphi dc\le I(\varphi)$.

Suppose $b=\int_X^{\vee+} \varphi dc<I(\varphi)=a$. Put $\varepsilon=(a-b)/2$.  Then for each $t\in\varphi(X)$ there exists $\chi^t\in\Upsilon^t_{\varphi_t}$ such that $I(\chi^t)<a-\epsilon$.  We can assume that $\chi^t$ is comonotone with $\varphi$ and $\chi^t$ refines $\varphi$ by Lemma \ref{Comon1}.  The set $V_t=\{y\mid \chi^t(y)+\varepsilon>\varphi(y)\}$ is an open neighborhood for each $x\in X$ with $\varphi(x)=t$. We can choose a finite subcover $\{V_{t_1},\dots,V_{t_k}\}$ of the open cover $\{V_t\mid t\in\varphi(X)\}$ of $X$. Since $\varphi$ refines each $\chi^{t_i}$, the functions $\chi^{t_i}+\varepsilon$ are pairwise comonotone, hence $I(\bigvee_{i=1}^k(\chi^{t_i}+\varepsilon))<a$. On the other hand $\varphi\le \bigvee_{i=1}^k(\chi^{t_i}+\varepsilon)$ and we obtain a contradiction.

Let us show uniqueness of the capacity $c$. Consider any capacity  $k\in MX$ such that $I(\varphi)=\int_X^{\vee+} \varphi dk$ for each $\varphi\in C(X)$. Suppose there exists $A\in\F(X)$ such that $k(A)\ne c(A)$. We can suppose that $k(A)<c(A)$. The proof is analogous in the opposite case. 

Take any $s\in\R$ such that $k(A)<s<c(A)$. Since the capacity $k$ is upper semicontinuous, there exists an open set $V\subset X$ such that $A\subset V$ and we have $k(B)<s$ for each compactum $B$ with $B\subset V$.
Consider a function $\psi\in C(X,[\ln s,0])$ such that $\psi|X\setminus V\equiv s$ and  $\psi|A\equiv 0$. Then we have $$\int_X^{\vee+} \varphi dc\ge \ln(c(\psi_0))+0\ge\ln(c(A)).$$ On the other hand we have $$\ln(k(\psi_t))+t\le \ln s<\ln(c(A))$$ for each $t\in\R$. We obtained a contradiction. 
\end{proof}

It is worth noting that the question whether we can eliminate monotonicity from the above characterization is a part of more general problem stated in \cite{Ra6}: if a  comonotonically maxitive functional is  monotone? This problem was affirmatively solved for finite compacta in  \cite{Ra6}.

\section{Categorical aspects.} We denote by $SX$ the set of all normalized, monotone,  comonotonically maximative and plus-homogeneous functionals on $C(X)$ for a compactum $X$. We consider $SX$ as a subspace of $\R^{C(X)}$. 

By $\Comp$ we denote the category of compact Hausdorff
spaces (compacta) and continuous maps. We recall the notion  of monad (or triple) in the sense of S.Eilenberg and J.Moore \cite{EilMoore}.  We define it only for the category $\Comp$.

A {\it monad} $\E=(E,\eta,\mu)$ in the category $\Comp$ consists of an endofunctor $E:{\Comp}\to{\Comp}$ and natural transformations $\eta:\Id_{\Comp}\to E$ (unity), $\mu:E^2\to E$ (multiplication) with components $\eta X:X\to EX$ and $\mu X:E^2X\to EX$ for which the following two diagrams are commutative for each compactum $X$
$$ \xymatrix{
  EX \ar[dr]_{\id_{EX}} \ar[r]^{E(\eta X)}
                & E^2X \ar[d]^{\mu X}&EX\ar[l]_{E(\eta X)}\ar[dl]^{\id_{EX}}  \\
                & EX             } $$
             and
  $$\xymatrix{
    E^3X \ar[d]_{\mu(EX)} \ar[r]^{E(\mu X)} & E^2X \ar[d]^{\mu X} \\
    E^2X\ar[r]^{\mu X} & EX  }$$

  (By $\Id_{\Comp}$ we denote the identity functor on the category ${\Comp}$ and $E^2$ is the superposition $E\circ E$ of $E$.)

We describe below some known monads in the category $\Comp$ based on constructions $O$, $I$ and $\Pi$ introduced in Section \ref{Pr}.  For a function $\phi\in C(X)$ we denote by $\pi_\phi$ or $\pi(\phi)$ the
corresponding projection $\pi_\phi:OX\to \R$. The construction $O$ is  functorial what means that for each continuous map $f:X\to Y$ we can consider a continuous map $Of:OX\to OY$ defined as follows $\pi_\psi\circ Of=\pi(\psi\circ f)$ for  $\psi\in C(Y)$.
The functor $O$ was completed to the monad $\OM=(O,\eta,\mu)$  in \cite{Ra1}. Let us describe the components of the natural transformations $\eta$ and $\mu$.  Let us remark that $\pi_\phi\in C(OX)$ and we can consider the map $\pi(\pi_\phi):O(OX)=O^2X\to\R$.  For a compactum $X$ we define components $\eta X$ and $\mu X$ of natural transformations $\eta:\Id_{\Comp}\to O$, $\mu:O^2\to O$ by $\pi_\phi\circ \eta X=\phi$ and $\pi_\phi\circ \mu X=\pi(\pi_\phi)$ for all $\phi\in C(X)$. It was proved in \cite{Ra1} that the triple $\OM=(O,\eta,\mu)$ forms a monad in the category $\Comp$.

It was shown in \cite{Za1} that the space $IX$ is compact and Hausdorff for a compactum $X$ and $Of(IX)\subset IY$ for each continuous map $f:X\to Y$, thus $I$ is a subfunctor of the functor $O$.  Moreover we have $\eta X(X)\subset IX$ and $\mu X(I^2X)\subset IX$ for each compactum $X$, hence the triple $\I=(I,\eta,\mu')$ is a submonad of the monad $\OM$, where $\mu':I^2\to I$ is the natural transformation with components $\mu' X=\mu X|I^2 X$.

The construction $\Pi$ was completed to the monad $\U_\bul=(\Pi,\eta,\mu_\bul)$ (where $\bul$ is the usual multiplication operation) in \cite{NR}.  For a continuous map of compacta $f:X\to Y$ we define the map $f:\Pi X\to \Pi Y$ by the formula $\Pi f(\nu)(A)=\nu(f^{-1}(A))$ where $\nu\in \Pi X$ and $A$ is a closed subset of $Y$. The map $\Pi f$ is continuous.  In fact, this extension of the construction $\Pi$ defines the possibility capacity functor $\Pi$ in the category $\Comp$.

The components of the  natural transformations $\eta$ and $\mu_\bul$ are defined as follows:
$$
\eta X(x)(F)=\begin{cases}
1,&x\in F,\\
0,&x\notin F;\end{cases}
$$

For a closed set $F\subset X$ and for $t\in [0,1]$ put $F_t=\{c\in MX\mid c(F)\ge t\}$. Define the map $\mu_\bul X:\Pi^2 X\to \Pi X$  by the formula $$\mu_\bul X(\C)(F)=\max\{\C(F_t)\cdot t\mid t\in(0,1]\}$$ for a closed set $F\subset X$ and $\C\in\Pi^2 X$. Let us remark that it is impossible to extend the monad structure $\U_\bul$ to the whole capacity functor $M$ \cite{Ra8}. However, $M$ can be completed to another monad based on the maximum  and the minimum operations \cite{NZ}.

An isomorphism of monads $\U_\bul$ and $\I$  was built in \cite{Ra6} and \cite{Ra5}. The main goal of this section is to extend this isomorphism to some correspondence between constructions $M$ and $S$ and investigate its categorical properties.  

For $\psi\in C(X)$ we define the function $lX^\psi:MX\to\R$ by the formula $lX^\psi(\nu)=\int_X^{\vee+} \psi d\nu$ for $\nu\in MX$. 

\begin{lemma}\label{contint} The map $lX^\psi$ is continuous for each $\psi\in C(X)$.
\end{lemma}

\begin{proof} Consider any $\nu\in MX$ such that $lX^\psi(\nu)<a$ for some $a\in\R$. Choose $s,p\in \R$ such that $\psi(X)\subset [s,p]$. Put $\varepsilon=a-lX^\psi(\nu)$.  Choose $k\in\N$ such that $\frac{p-s}{k}<\varepsilon/2$ and put $t_i=s+\frac{i(p-s)}{k}$ for $i\in\{0,\dots,k\}$. For each $i\in\{0,\dots,k\}$ choose $\delta_i>0$ such that $\ln(\nu(\psi_{t_i})+\delta_i)<\ln(\nu(\psi_{t_i}))+\varepsilon/2$.  Define an open set $O_i=\{\mu\in MX\mid \mu(\psi_{t_i})<\nu(\psi_{t_i})+\delta_i\}$ and put $O=\cap_{i=1}^kO_i$. Evidently $O$ is an open neighborhood of $\nu$. Consider any $\mu\in O$ and $t\in[s,p]$. Let $i$ be a maximal element of $\{0,\dots,k\}$ such that $t_i\le t$. Then we have $\ln(\mu(\psi_{t}))+ t<\ln(\mu(\psi_{t_i}))+ t_i+\varepsilon/2<\ln(\nu(\psi_{t_i})+\delta_i)+ t_i+\varepsilon/2<\ln(\nu(\psi_{t_i}))+ t_i+\varepsilon\le a$.
Hence $lX^\psi(\mu)=\max\{\ln(\mu(\varphi_t))+ t\mid t\in\R\}=\max\{\ln(\mu(\varphi_t))+ t\mid t\in[s,p]\}<a$.

Now, consider any $\nu\in MX$ such that $lX^\psi(\nu)>a$ for some $a\in\R$. Then there exists $t\in\R$ such that $\ln(\nu(\psi_{t}))+ t> a$. Put $\varepsilon=\ln(\nu(\psi_{t}))+ t- a$. We choose $\delta>0$ such that $\ln(\nu(\psi_{t})-\delta)>\ln(\nu(\psi_{t}))-\varepsilon/2$. We can suppose $\delta<\varepsilon/2$.

Define an open set $O=\{\mu\in MX\mid \mu(\psi^{-1}((t-\delta,+\infty)))>\nu(\psi_t)-\delta\}$. Evidently $O$ is an open neighborhood of $\nu$. Consider any $\mu\in O$. There exists $p\in(t-\delta,t]$ such that $\mu(\psi_p)>\nu(\psi^{-1}((t-\delta,+\infty)))-\delta\ge\nu(\psi_{t})-\delta$. Then we have $\ln(\mu(\psi_{p}))+ p> \ln(\nu(\psi_{t})-\delta)+ t-\varepsilon/2>\ln(\nu(\psi_{t}))-\varepsilon/2+ t-\varepsilon/2=a$. Hence $lX^\psi(\mu)>a$ and the map $lX^\psi$ is continuous.
\end{proof}

We  define the map $lX:MX\to\R^{C(X)}$ taking the diagonal product $lX=(lX^\psi)_{\psi\in C(X)}$. Lemma \ref{prop} implies that $lX(MX)\subset SX$. Theorem \ref{repr} yields that the map $lX:MX\to SX$ is bijective. Finally, the continuity of $lX$ follows immediately from Lemma \ref{contint}. Hence we have the following statement.

\begin{theorem}\label{homeo} The map $lX:MX\to SX$ is a homeomorphism.
\end{theorem}

Theorem \ref{homeo} implies compactness of $SX$. Let $f:X\to Y$ be a continuous map. It is easy to check that $Of(SX)\subset SY$, thus we can define the  map $Sf:SX\to SY$ as the restriction of $Of$ on $SX$. Hence we obtain the functor of  comonotonically maxitive functionals $S$. Let us remark that the idempotent measure functor $I$ is a subfunctor of $S$ and $S$ is a subfunctor of the functor of order-preserving functionals $O$. 

It is easy to check that the maps $lX$ are components of the natural transformation $l:M\to S$ and we obtain that the functors $M$ and $S$ are isomorphic. Since the functors $I$ and $O$ are completed to the monads $\I$ and $\OM$ such that $\I$ is a submonad $\OM$, the question arises naturally, if the functor $S$ can be completed to a submonad of $\OM=(O,\eta,\mu)$. Let us remark that the affirmative answer to this question follows from the problem if we have the inclusion $\mu_X(S^2X)\subset S(X)$.

\end{document}